\documentclass[11pt,a4paper]{article}

\usepackage[hmargin=24mm,vmargin=22mm]{geometry}
\usepackage{newtxtext}
\usepackage{amsmath,amsthm}
\usepackage{newtxmath}
\usepackage{microtype}

\usepackage{xcolor}
\usepackage[unicode,colorlinks=true,allcolors=blue!45!black]{hyperref}

\newtheorem{theorem}{Theorem}[section]
\newtheorem{lemma}[theorem]{Lemma}
\newtheorem{corollary}[theorem]{Corollary}
\newtheorem{proposition}[theorem]{Proposition}
\theoremstyle{definition}
\newtheorem{definition}[theorem]{Definition}
\newtheorem{example}[theorem]{Example}

\newcommand{\Root}[1]{\sqrt{#1}}
\newcommand{\Roots}[2]{\mathcal R_{#1}(#2)}
\DeclareMathOperator{\FP}{FP}
\DeclareMathOperator{\FC}{FC}

\title{Dense Exact Factorizations\\
      of Countable Topological Groups}
\author{Mikhail Kabenyuk}
\date{}

\begin{document}
\maketitle
\begin{abstract}
We prove that every countably infinite group endowed with a topology
of countable weight in which every nonempty open set is infinite
admits an exact factorization into two dense infinite subsets.
More generally, every finite partial exact factorization with
nonempty factors extends to an exact factorization of the whole
group in which both factors meet every member of a prescribed
at most countable family of infinite subsets in infinitely many points. 
This settles the countable case of a question of Protasov
and Slobodianiuk. 
The proof uses the finite graph Ramsey theorem,
invariant finitely additive probability measures on amenable
groups, and a structural lemma on quadratic coverings based
on Neumann's coset covering lemma.
\end{abstract}



\section*{Introduction}

Exact factorizations of groups into products of subsets are studied
in connection with unique representations of elements, tilings,
and combinatorial properties of groups; a general account of this
subject can be found in the monographs~\cite{Szabo,SzaboSands2009}.
In the topological version of the problem, uniqueness of representation
is supplemented by density requirements on the factors.

Density of both factors was studied by Protasov and
Slobodianiuk~\cite{PS2017}. They proved that if a group $G$ has
uncountable cardinality $\kappa$ and is endowed with a topology
of weight at most $\kappa$ in which every nonempty open set has
cardinality $\kappa$, then there is an exact factorization $G=AB$
with dense factors of cardinality $\kappa$.
The countable case was left open.
In particular, the authors asked whether such a factorization exists
for every countable nondiscrete Hausdorff left topological group
of countable weight. For countable topological groups,
positive results were obtained in~\cite{PS2017} for abelian groups,
for groups in which all finitely generated subgroups are nowhere dense, 
and for groups in which the set of squares of elements of
every nonempty open set is infinite.

Another type of restriction on the factors was investigated by
Protasov~\cite{Protasov2021}. He proved that every countably infinite
topologizable group, that is, a group admitting a nondiscrete Hausdorff
group topology, has an exact factorization into two thin subsets.
A subset $X\subseteq G$ is thin if $gX\cap X$ and $Xg\cap X$ are
finite for every $g\ne e$. The same paper asks whether every
countably infinite group admits an exact factorization into
two infinite subsets. Protasov also posed this question in the
Kourovka Notebook~\cite[Problem~20.73]{Kourovka2026}.

Di~Siena~\cite{DiSiena2026} gave an affirmative answer to this algebraic question. 
The proof of Lemma~3.1 of~\cite{DiSiena2026} extends finite partial factorizations under
the assumption that every nonidentity element has only finitely
many square roots.
Here we remove the finiteness assumption from this extension step:
Lemma~\ref{lem:extension} applies to every infinite group.

We answer the countable question of~\cite{PS2017} in its general form: 
every countably infinite group endowed with a topology
having a countable base in which all nonempty open sets are infinite
admits an exact factorization into two dense infinite subsets.
Theorem~\ref{thm:main} proves more: 
any finite partial exact
factorization extends to an exact factorization whose two factors
each meet every member of a prescribed at most countable family
of infinite subsets in infinitely many points.

The key step is condition (R), concerning finite unions of translates 
of fibers of the squaring map. 
We first prove it for amenable groups using the finite graph
Ramsey theorem and a finitely additive probability measure
invariant on both sides.
Independently, Lemma~\ref{lem:quadratic-cover}, based on Neumann's
lemma on coverings by cosets, shows that a finite quadratic covering
outside a finite set forces the FC-center to have finite index.
If (R) failed, this structural conclusion would make the group
amenable, contradicting the previously established estimate.
This proves (R) for all infinite groups and supplies the extension
step needed for the countable construction.

\section{Main result and notation}

We denote the identity element of a group by $e$ and the set of
positive integers by $\mathbb{N}$.
A product $AB$ of subsets of a group $G$ is called a
\emph{partial exact factorization} if the map
$A\times B\longrightarrow G$, $(a,b)\longmapsto ab$, is injective.
If, in addition, $AB=G$, we call it a \emph{complete exact
factorization} of $G$, or simply an \emph{exact factorization of $G$}.
In this case, every element of $G$ has exactly one
representation as $ab$, where $a\in A$ and $b\in B$.

\begin{theorem}\label{thm:main}
Let $G$ be a countably infinite group, and let
$\mathcal X$ be an at most countable family of infinite subsets of $G$. 
Every partial exact factorization $A_0B_0$ with finite
nonempty $A_0,B_0$ extends to an exact factorization
\[
 G=AB,\quad A_0\subseteq A,\quad B_0\subseteq B,
\]
where $A,B$ are countably infinite and
\[
 |A\cap X|=|B\cap X|=\aleph_0\quad(\forall X\in\mathcal X).
\]
\end{theorem}

For $h\in G$ and $S\subseteq G$, put
\[
 \Roots{G}{h}=\{x\in G:x^2=h\},\quad
 \Root{S}=\{x\in G:x^2\in S\}=\bigcup_{h\in S}\Roots{G}{h}.
\]
We call $\Roots{G}{h}$ the \emph{fiber of the squaring map over $h$}.
Fibers over distinct elements are disjoint.
For an infinite sequence $(x_n)_{n\geq1}$ in a group, write
\[
 \FP(x_n)=\{x_{i_1}\cdots x_{i_k}:k\geq1,\ i_1<\cdots<i_k\}.
\]
This is the set of all nonempty finite products with increasing indices. 
The terms need not be distinct. This notation will also
be used when the ambient group is finite.

\begin{samepage}
For the proof, we introduce the following auxiliary condition.

\begin{definition}\label{def:R}
An infinite group $G$ satisfies condition \textup{(R)} if, for
every choice of pairwise distinct elements $h_1,\ldots,h_r\ne e$,
where $r\geq1$, and arbitrary $t_1,\ldots,t_r\in G$, the set
\[
 G\setminus\bigcup_{i=1}^r\Roots{G}{h_i}t_i
\]
is infinite.
\end{definition}
\end{samepage}
We first prove condition (R) for amenable groups, and then for all
infinite groups. Condition (R) will then allow us to extend finite
partial factorizations and carry out a countable induction.

The following examples illustrate condition~\textup{(R)} and the
roles of its restrictions on the elements $h_i$.

\begin{example}[An infinite direct sum of groups of order two]
\label{ex:R-boolean}
Let
\[
 E=\bigoplus_{n=1}^{\infty}C_2.
\]
Every element of $E$ has square $e$. Thus
$\Roots{E}{h}=\varnothing$ for every $h\ne e$, and the set in
Definition~\ref{def:R} is the whole group:
\[
 E\setminus\bigcup_{i=1}^r\Roots{E}{h_i}t_i=E.
\]
This example also explains why the condition $h_i\ne e$ is
necessary: $\Roots{E}{e}t=E$ for every $t\in E$.
\end{example}

\begin{example}[The infinite dihedral group]
\label{ex:R-dihedral}
Let
\[
 D_\infty=\langle a,s\mid s^2=e,\ sas=a^{-1}\rangle.
\]
Since $a$ has infinite order, for $h\ne e$ we have
\[
 \Roots{D_\infty}{h}=
 \begin{cases}
  \{a^k\},& h=a^{2k},\quad k\in\mathbb Z\setminus\{0\},\\
  \varnothing,&\text{otherwise}.
 \end{cases}
\]
Hence the excluded union contains at most $r$ elements, for instance,
\[
 D_\infty\setminus\Roots{D_\infty}{a^2}s
   =D_\infty\setminus\{as\}.
\]
\end{example}

\begin{example}[An infinite fiber over a nonidentity element]
\label{ex:R-C4-boolean}
Let $E$ be as in Example~\ref{ex:R-boolean}, and put
\[
 G=C_4\times E,\quad C_4=\langle a\rangle,\quad
 H=\langle a^2\rangle\times E.
\]
We identify the factors with their natural subgroups of $G$.
Then $H$ is an infinite subgroup of index two,
\[
 \Roots{G}{e}=H,\quad \Roots{G}{a^2}=aH,
\]
and all other fibers are empty. Since the elements $h_i$ are
pairwise distinct and different from $e$, at most one of the
fibers $\Roots{G}{h_i}$ is nonempty. The excluded union is
therefore either empty or a single coset of $H$.
Its complement is either $G$ or the other coset of $H$,
and is infinite in both cases. In particular,
\[
 G\setminus\Roots{G}{a^2}=H,
\]
so here both the excluded set and its complement are infinite.

The requirement that the elements $h_i$ be pairwise distinct
is essential. If repetitions were allowed, two translates of
the same fiber would cover the group:
\[
 r=2,\ h_1=h_2=a^2,\ t_1=e,\ t_2=a\quad\Rightarrow\quad
 \Roots{G}{a^2}\cup\Roots{G}{a^2}a=aH\cup H=G.
\]
\end{example}

\section{A combinatorial lemma on square roots}\label{sec:ramsey}

\subsection{The finite graph Ramsey theorem}

We use the following standard form of Ramsey's
theorem~\cite{Ramsey1930}. An edge coloring with $r$ colors assigns
one of $r$ values to each edge. A subset of the vertices of a graph
is called a clique if each pair of its vertices is joined by an edge.
A clique is monochromatic if all its edges have the same color.

\begin{theorem}[Finite Ramsey theorem]\label{thm:ramsey}
For any positive integers $r$ and $m$, there exists a positive
integer $N$ such that every complete graph on $N$ vertices
whose edges are colored with $r$ colors contains
a monochromatic clique on $m$ vertices.
\end{theorem}

\subsection{An algebraic observation}

\begin{lemma}\label{lem:commute}
Let $u,v,w$ be elements of a group such that
\[
 u^2=v^2=w^2=a,\quad (uw)^2=(vw)^2.
\]
Then $uv$ commutes with $w$.
\end{lemma}

\begin{proof}
The element $a$ commutes with $u,v,w$, since it is the square of
each of them. Cancelling $w$ on the right in $(uw)^2=(vw)^2$,
we obtain $uwu=vwv$. Multiplying on the left by $u$ and on the
right by $v$ gives
\[
 a\,wuv=uvw\,a=a\,uvw.
\]
Cancelling $a$ yields $wuv=uvw$.
\end{proof}

\subsection{A lemma on finite products}

\begin{lemma}\label{lem:no-ip}
Let $G$ be an arbitrary group, and let
$S\subseteq G\setminus\{e\}$ be finite. Then $\Root{S}$ does not
contain $\FP(x_n)$ for any infinite sequence $(x_n)$.
\end{lemma}

\begin{proof}
The assertion is clear when $S=\varnothing$. Let $r=|S|\geq1$,
and suppose that $\FP(x_n)\subseteq\Root{S}$.
Thus the square of every nonempty finite product with increasing
indices belongs to $S$.

If the sequence $(x_n)$ takes only finitely many values, some
element $g$ occurs infinitely often. Then $g^k\in\FP(x_n)$ and
$g^{2k}\in S$ for every $k\geq1$. Since $S$ is finite,
$g^{2p}=g^{2q}$ for some $p<q$, whence
$e=g^{2(q-p)}\in S$, a contradiction. We may therefore pass to
a subsequence $(x_n)$ whose terms are pairwise distinct.
Passing to a further subsequence, we may also assume that
\[
 x_n^2=a\quad(n\geq1)
\]
for some $a\in S$.

In Theorem~\ref{thm:ramsey}, take $m=2r+2$ and the corresponding
value of $N$. Form the complete graph with vertices
$x_1,\ldots,x_N\in G$. Assign to the edge $\{x_i,x_j\}$, where
$i<j$, the color
\[
 c(\{x_i,x_j\})=(x_ix_j)^2\in S.
\]
The graph is undirected; the order of the factors is determined by
the order of the elements in the original sequence.
A monochromatic clique gives elements $y_1,\ldots,y_{2r+2}$,
listed in the order of their original indices, such that
\[
 y_i^2=a,\quad (y_iy_j)^2=b\quad(i<j)
\]
for some $b\in S$. The values $a$ and $b$ need not coincide.
By Lemma~\ref{lem:commute}, the product $y_iy_j$ commutes with
$y_k$ whenever $i<j<k$.

Partition the chosen $2r+2$ elements into $r+1$ consecutive pairs:
\[
 z_\ell=y_{2\ell-1}y_{2\ell}\quad(1\leq\ell\leq r+1).
\]
Each $z_\ell$ has square $b$. If $\ell<t$, then $z_\ell$ commutes
with both $y_{2t-1}$ and $y_{2t}$, and hence with $z_t$.
Consequently, $z_1,\ldots,z_{r+1}$ commute pairwise and
\[
 (z_1\cdots z_k)^2=b^k\quad(1\leq k\leq r+1).
\]
But $z_1\cdots z_k=y_1\cdots y_{2k}\in\FP(x_n)$.
Therefore $b,b^2,\ldots,b^{r+1}\in S$.
This is precisely where the $r+1$ pairs are needed: at least two
of these $r+1$ powers coincide, since $|S|=r$.
Let $b^p=b^q$, where $1\leq p<q\leq r+1$.
Then $b^{q-p}=e$, and this power also belongs to $S$, since
$1\leq q-p\leq r$. We obtain the contradiction $e\in S$.
\end{proof}

\section{Amenability and finitely additive measures}\label{sec:amenable}

\subsection{Definitions and standard facts}

Throughout, amenability refers to the underlying discrete group.
The standard facts used below can be found in Juschenko's
monograph~\cite{Juschenko2022}.

\begin{definition}
A \emph{finitely additive probability measure} on a group $G$ is
a function $\mu:\mathcal P(G)\to[0,1]$, defined on all subsets
of $G$, such that
\[
 \mu(G)=1,\quad
 \mu(E\cup F)=\mu(E)+\mu(F)\quad(E\cap F=\varnothing).
\]
The group $G$ is \emph{amenable} if it admits such a measure
that is \emph{left invariant}, meaning that
\[
 \mu(gE)=\mu(E)\quad(\forall g\in G,\ E\subseteq G).
\]
\end{definition}
In an amenable group, the measure can be chosen to be
\emph{invariant on both sides}:
\[
 \mu(gE)=\mu(Eg)=\mu(E)\quad(\forall g\in G,\ E\subseteq G).
\]
We shall use the following consequences of finite additivity:
\begin{gather*}
 \mu(G\setminus E)=1-\mu(E),\qquad
 \mu\left(\bigcup_{i=1}^kE_i\right)\leq\sum_{i=1}^k\mu(E_i),\\
 \mu(E_1)=\cdots=\mu(E_k)=1
 \quad\Longrightarrow\quad
 \mu\left(\bigcap_{i=1}^kE_i\right)=1.
\end{gather*}
On an infinite group, every left invariant finitely additive
probability measure vanishes on finite sets. Hence every set of
positive measure is infinite, and removing a finite set does not
change its measure.

Finite and abelian groups are amenable.
The class of amenable groups is closed under taking subgroups,
homomorphic images, extensions, and directed unions of subgroups.
Here closure under extensions means that if
$N\mathrel{\trianglelefteq}G$ and both $N$ and $G/N$ are amenable,
then $G$ is amenable. 
It follows that finite direct products of amenable groups are amenable, 
as are groups containing an amenable subgroup of finite index.

\subsection{The measure of a set of square roots}

The following lemma uses the combinatorial result of
Section~\ref{sec:ramsey}.

\begin{lemma}\label{lem:measure}
Let $G$ be an amenable group, finite or infinite, let $\mu$ be a
left invariant finitely additive probability measure on all its
subsets, and let $S\subseteq G\setminus\{e\}$ be finite. Then
\[
 \mu(\Root{S})<1.
\]
\end{lemma}

\begin{proof}
Suppose that $\mu(T)=1$, where $T=\Root{S}$.
Each translate $p^{-1}T$ has measure $1$, as does every finite
intersection of such translates. In particular, these intersections
are nonempty.

Choose $x_1\in T$. If $x_1,\ldots,x_n$ have already been chosen
and all their nonempty products with increasing indices belong
to $T$, put
\[
 K_n=\{e\}\cup
 \{x_{i_1}\cdots x_{i_k}:1\leq k\leq n,
                      \ 1\leq i_1<\cdots<i_k\leq n\}
\]
and choose
\[
 x_{n+1}\in\bigcap_{p\in K_n}p^{-1}T.
\]
For every $p\in K_n$, we have $px_{n+1}\in T$.
Taking $p=e$ gives $x_{n+1}\in T$, while the remaining values
of $p$ give all new products ending in $x_{n+1}$.
The old products remain in $T$ by the induction hypothesis.
This yields an infinite sequence with $\FP(x_n)\subseteq T$,
contrary to Lemma~\ref{lem:no-ip}.
\end{proof}

\subsection{Verification of condition (R)}

\begin{proposition}\label{prop:amenable}
Let $G$ be amenable, let $h_1,\ldots,h_r\ne e$ be pairwise
distinct, and let $t_1,\ldots,t_r\in G$. For every finitely additive
probability measure $\mu$ on $G$ that is invariant on both sides,
we have
\[
 \mu\left(\bigcup_{i=1}^r\Roots{G}{h_i}t_i\right)<1.
\]
In particular, this union is not equal to $G$.
If $G$ is an infinite amenable group, then it satisfies \textup{(R)}.
\end{proposition}

\begin{proof}
Put $S=\{h_1,\ldots,h_r\}$.
By the properties of $\mu$ and Lemma~\ref{lem:measure},
\[
 \mu\left(\bigcup_{i=1}^r\Roots{G}{h_i}t_i\right)
 \leq\sum_{i=1}^r\mu(\Roots{G}{h_i})
 =\mu(\Root{S})<1.
\]
The inequality follows from subadditivity and right invariance,
and the equality follows from the pairwise disjointness of the
untranslated fibers. The translated fibers may intersect.
The complement has positive measure; if $G$ is infinite,
the complement is infinite. This is precisely condition (R).
\end{proof}

\section{Quadratic coverings and the general case of (R)}

Recall the classical lemma of
B.~H.~Neumann~\cite[Lemma~4.1]{Neumann1954}:
if a group is the union of finitely many cosets of subgroups,
then at least one of these subgroups has finite index.
Left and right cosets can be treated alike, since
$Hg=g(g^{-1}Hg)$.

For a group $G$, let
\[
 \Delta=\FC(G)=\{c\in G:|G:C_G(c)|<\infty\}
\]
denote its \emph{FC-center}, the subgroup of elements with finite
conjugacy classes. This subgroup is normal. Indeed, if the
centralizers of $c$ and $d$ have finite index, their intersection
has finite index and is contained in the centralizer of $cd$.
The corresponding conditions are also preserved under taking
inverses and conjugates.

\begin{samepage}
\begin{lemma}\label{lem:fc-amenable}
The FC-center of every group is amenable. If it has finite
index, then the whole group is amenable.
\end{lemma}

\begin{proof}
Let $\Delta=\FC(G)$. Every finitely generated subgroup
$H=\langle u_1,\ldots,u_m\rangle\leq\Delta$ has a center
of finite index. Indeed,
\[
 Z(H)=\bigcap_{k=1}^m\bigl(H\cap C_G(u_k)\bigr),
\]
and each subgroup $H\cap C_G(u_k)$ has finite index in $H$,
since $u_k\in\FC(G)$.
Thus $H$ is amenable, being a group with an abelian subgroup
of finite index. As a directed union of such subgroups,
$\Delta$ is amenable.

If $|G:\Delta|<\infty$, then $G$ is amenable, since
it contains an amenable subgroup of finite index.
\end{proof}
\end{samepage}

We now establish a structural consequence of a quadratic covering.

\begin{lemma}\label{lem:quadratic-cover}
Let $G$ be infinite, let $F\subseteq G$ be finite, and let
$a_1,\ldots,a_r,b_1,\ldots,b_r\in G$, where $r\geq1$.
Suppose that for every $x\in G\setminus F$ there is at least one
$i\in\{1,\ldots,r\}$ such that
\[
 xa_ix=b_i.
\]
Then $|G:\FC(G)|\leq r^2$.
\end{lemma}

\begin{proof}
Fix an arbitrary $g\in G$.
If $x\notin F\cup g^{-1}F$, then both $x$ and $gx$ lie outside
$F$. Hence there are indices $i,j$ such that
\begin{equation}\label{eq:pair}
 xa_ix=b_i,\quad (gx)a_j(gx)=b_j.
\end{equation}
The first equality gives
\[
 x=a_i^{-1}x^{-1}b_i,
\]
and the second gives $xa_jgx=g^{-1}b_j$.
Substituting this expression for the last $x$, we obtain
\[
 xa_jga_i^{-1}x^{-1}b_i=g^{-1}b_j,
\]
that is,
\begin{equation}\label{eq:conjugator}
 x\bigl(a_jga_i^{-1}\bigr)x^{-1}=g^{-1}b_jb_i^{-1}.
\end{equation}
Thus, for fixed $g$, the set $G\setminus(F\cup g^{-1}F)$ is covered
by at most $r^2$ solution sets of equations~\eqref{eq:conjugator}.

If nonempty, the solution set of $xcx^{-1}=d$ is a left coset
$x_0C_G(c)$, where $x_0$ is any one solution. Indeed, $x$ is a
solution if and only if $x_0^{-1}x$ commutes with $c$.
Adding the singleton sets corresponding to elements of
$F\cup g^{-1}F$, we obtain a covering of the whole group $G$
by finitely many cosets of subgroups. Singleton sets are cosets
of the trivial subgroup, which has infinite index since $G$
is infinite. By Neumann's lemma, at least one centralizer
\[
 C_G(a_jga_i^{-1})
\]
has finite index. Consequently, for some $i,j$,
\[
 a_jga_i^{-1}\in\Delta,
 \quad
 g\in a_j^{-1}\Delta a_i=a_j^{-1}a_i\Delta.
\]
Here we use the normality of $\Delta$.
Since $g$ was arbitrary, at most $r^2$ cosets
$a_j^{-1}a_i\Delta$ cover $G$.
This proves $|G:\Delta|\leq r^2$.
\end{proof}

\begin{samepage}
\begin{theorem}\label{thm:R-all}
Every infinite group satisfies \textup{(R)}.
\end{theorem}

\begin{proof}
Suppose, to the contrary, that for some pairwise distinct
$h_1,\ldots,h_r\ne e$
and some $t_1,\ldots,t_r\in G$, the complement of
\[
 D=\bigcup_{i=1}^r\Roots{G}{h_i}t_i
\]
is finite. For every $x\in D$, there is an $i$ such that
\[
 (xt_i^{-1})^2=h_i
 \quad\Longleftrightarrow\quad
 xt_i^{-1}x=h_it_i.
\]
By Lemma~\ref{lem:quadratic-cover}, applied with $a_i=t_i^{-1}$
and $b_i=h_it_i$, the subgroup $\Delta=\FC(G)$ has finite index in $G$. 
By Lemma~\ref{lem:fc-amenable}, the group $G$ is amenable.
Since $G$ is infinite, Proposition~\ref{prop:amenable}
implies that $G$ satisfies \textup{(R)}, contrary to our
assumption.
\end{proof}
\end{samepage}

\section{Construction of the factorization}

\subsection{Extending a partial factorization}

An inductive construction using finite partial factorizations
was employed in~\cite{PS2017}. The device of simultaneously
adding $x$ and $x^{-1}g$ to the two factors appears in the proof
of Lemma~3.1 of~\cite{DiSiena2026}, under the assumption that
all fibers of the squaring map over nonidentity elements are finite.
The following lemma justifies this step for an arbitrary infinite
group. The reduction of possible collisions to quadratic equations
is retained, while the existence of an admissible value of $x$
is ensured by Theorem~\ref{thm:R-all}, which does not require
the fibers to be finite.

\begin{lemma}\label{lem:extension}
Let $G$ be an infinite group, let $AB$ be a partial exact
factorization with finite nonempty $A,B$, and let $g\notin AB$.
Then there exists $x\in G$ such that the product
\[
 \bigl(A\cup\{x\}\bigr)\bigl(B\cup\{x^{-1}g\}\bigr)
\]
is also a partial exact factorization and contains $g$.
The element $x$ can be chosen outside any prescribed finite set.
\end{lemma}

\begin{proof}
By Theorem~\ref{thm:R-all}, the group $G$ satisfies (R).
Put $y=x^{-1}g$. After the extension, all products fall into
four parts:
\begin{equation}\label{eq:four}
 AB,\quad xB,\quad Ay,\quad\{g\}.
\end{equation}
Choosing $x$ outside the finite set
\[
 E=ABB^{-1}\cup gB^{-1}A^{-1}A
\]
ensures that
\[
 x\notin A,\quad y\notin B,\quad
 xB\cap AB=\varnothing,\quad Ay\cap AB=\varnothing.
\]
Indeed, the equalities $xb=a'b'$ and $ay=a'b'$ give,
respectively, $x=a'b'b^{-1}\in ABB^{-1}$ and
$x=g(b')^{-1}(a')^{-1}a\in gB^{-1}A^{-1}A$.
The inclusions $A\subseteq E$ and $gB^{-1}\subseteq E$ follow
from the fact that $A,B$ are nonempty. The equalities $g=xb$
and $g=ay$ would imply $y=b\in B$ and $x=a\in A$,
respectively, while $g\notin AB$ by assumption.

It remains to ensure that $xB\cap Ay=\varnothing$.
For $a\in A$, $b\in B$,
we have
\begin{equation}\label{eq:quadratic}
 xb=ax^{-1}g
 \quad\Longleftrightarrow\quad
 (xbg^{-1})^2=abg^{-1}.
\end{equation}
Thus the set of all values of $x$ forbidden by these collisions is
\[
 D=\bigcup_{\substack{a\in A\\b\in B}}
       \Roots{G}{abg^{-1}}\,gb^{-1}.
\]
The elements $abg^{-1}$ are pairwise distinct because $AB$ is
a partial exact factorization, and they are not equal to $e$
because $g\notin AB$. By (R), the set $G\setminus D$ is infinite.
Deleting $E$ and any further prescribed finite set still leaves
infinitely many choices for $x$. Choose one of them.
Then the parts in~\eqref{eq:four} are pairwise
disjoint, and within each part distinct pairs of factors give
distinct products. Hence the extended product is a partial
exact factorization.
\end{proof}

\subsection{Proof of the main theorem}

\begin{proof}[Proof of Theorem~\ref{thm:main}]
Enumerate $G=\{g_0,g_1,\ldots\}$.
Choose a sequence $(X_n)_{n\geq0}$ with values in
$\mathcal X\cup\{G\}$
in which every member of this family occurs infinitely often.
We construct increasing finite sets $A_n,B_n$, starting with
the prescribed $A_0,B_0$ and preserving the property that
each $A_nB_n$ is a partial exact factorization.

At stage $n$, first choose
\[
 a_n\in X_n\setminus A_nB_nB_n^{-1},\quad
 \widetilde A_n=A_n\cup\{a_n\}.
\]
The excluded set is finite, so this choice is possible.
An equality $a_nb=ab'$ with $a\in A_n$ and $b,b'\in B_n$ would imply
$a_n=ab'b^{-1}$, which has been excluded. Thus
$\widetilde A_nB_n$ is a partial exact factorization, and
$a_n\notin A_n$. Next, choose
\[
 b_n\in X_n\setminus\widetilde A_n^{-1}\widetilde A_nB_n,
 \quad \widetilde B_n=B_n\cup\{b_n\}.
\]
An equality $ab_n=a'b$ with $a,a'\in\widetilde A_n$ and
$b\in B_n$ would imply $b_n=a^{-1}a'b$, which has been excluded.
Thus $\widetilde A_n\widetilde B_n$ is a partial exact
factorization and $b_n\notin B_n$.

If $g_n$ already belongs to this product, put
$A_{n+1}=\widetilde A_n$, $B_{n+1}=\widetilde B_n$.
Otherwise, apply Lemma~\ref{lem:extension} and extend both
factors so that $g_n\in A_{n+1}B_{n+1}$.

Put $A=\bigcup_{n\geq0}A_n$, $B=\bigcup_{n\geq0}B_n$.
Every $g_n$ is represented, so $AB=G$.
If an element had two distinct representations in $AB$,
all four factors would occur at some finite stage, contradicting
the fact that $A_nB_n$ is a partial exact factorization.
Each $X\in\mathcal X\cup\{G\}$ is considered infinitely often,
and new elements are added to both factors each time.
This gives all the required infinite intersections.
\end{proof}

\subsection{Topological consequences}

\begin{corollary}\label{cor:topology}
Let $G$ be a countably infinite group endowed with a topology
having a countable base in which all nonempty open sets are infinite.
Then $G$ admits an exact factorization $G=AB$ into two dense
countably infinite subsets. Any prescribed partial exact
factorization with finite nonempty factors can be extended in
this way. Moreover, $A$ and $B$ can be chosen so that each has
infinite intersection with every nonempty open set.
\end{corollary}

\begin{proof}
In Theorem~\ref{thm:main}, take $\mathcal X$ to be the family of
nonempty members of a countable base. Every nonempty open set
contains one of these members, so all the asserted properties follow.
\end{proof}

\begin{corollary}\label{cor:left}
Every countable nondiscrete Hausdorff left topological group of
countable weight admits an exact factorization into two dense
countably infinite subsets.
\end{corollary}

\begin{proof}
In a left topological group, all left translations are continuous
and hence are homeomorphisms. If there were a nonempty finite
open set, the $T_1$ property would give an isolated point.
Left translations would then make every point isolated,
contrary to nondiscreteness. Thus all nonempty open sets are
infinite, and Corollary~\ref{cor:topology} applies.
\end{proof}

\noindent
This gives an affirmative answer to the question for countable groups
posed in~\cite{PS2017}.

\vspace{10mm}
\noindent\textbf{Use of AI.}
The author used OpenAI tools in developing the proofs
and preparing the manuscript.


\begin{thebibliography}{9}
\raggedright

\bibitem{DiSiena2026}
B.~G.~Di~Siena.
\emph{Exact factorizations of countable groups}.
Research preprint, ProofForum,
\href{https://www.proofforum.org/abs/PF:2026.000023/v1}
     {PF:2026.000023}, version~1, October 13, 2026.
Also available from the
\href{https://kourovkanotebookorg.wordpress.com/wp-content/uploads/2026/10/20.73.pdf}
     {Kourovka Notebook website}.

\bibitem{Juschenko2022}
K.~Juschenko.
\emph{Amenability of Discrete Groups by Examples}.
Mathematical Surveys and Monographs, vol.~266.
American Mathematical Society, Providence, RI, 2022.
\href{https://doi.org/10.1090/surv/266}{doi:10.1090/surv/266}.

\bibitem{Kourovka2026}
E.~I.~Khukhro, V.~D.~Mazurov (eds.).
\emph{Unsolved Problems in Group Theory. The Kourovka Notebook}.
No.~21. Sobolev Institute of Mathematics, Novosibirsk, 2026.
November update.
\href{https://kourovkanotebookorg.wordpress.com/wp-content/uploads/2026/11/21tkt.pdf}
     {Electronic version};
\href{https://kourovkanotebookorg.wordpress.com/}{updates}.

\bibitem{Neumann1954}
B.~H.~Neumann.
\emph{Groups covered by permutable subsets}.
Journal of the London Mathematical Society (1)
\textbf{29} (1954), no.~2, 236--248.
\href{https://doi.org/10.1112/jlms/s1-29.2.236}
     {doi:10.1112/jlms/s1-29.2.236}.

\bibitem{Protasov2021}
I.~Protasov.
\emph{On factoring groups into thin subsets}.
Axioms \textbf{10} (2021), no.~2, article~89.
\href{https://doi.org/10.3390/axioms10020089}
     {doi:10.3390/axioms10020089}.

\bibitem{PS2017}
I.~Protasov, S.~Slobodianiuk.
\emph{A note on factoring groups into dense subsets}.
Journal of Group Theory \textbf{20} (2017), no.~1, 33--38.
\href{https://doi.org/10.1515/jgth-2016-0021}
     {doi:10.1515/jgth-2016-0021}.
Preprint: \href{https://arxiv.org/abs/1602.01603}{arXiv:1602.01603}.

\bibitem{Ramsey1930}
F.~P.~Ramsey.
\emph{On a problem of formal logic}.
Proceedings of the London Mathematical Society (2)
\textbf{30} (1930), no.~1, 264--286.
\href{https://doi.org/10.1112/plms/s2-30.1.264}
     {doi:10.1112/plms/s2-30.1.264}.

\bibitem{Szabo}
S.~Szab\'o.
\emph{Topics in Factorization of Abelian Groups}.
Hindustan Book Agency, New Delhi, 2004.
\href{https://doi.org/10.1007/978-93-86279-22-4}
     {doi:10.1007/978-93-86279-22-4}.

\bibitem{SzaboSands2009}
S.~Szab\'o\ and A.~D.~Sands.
\emph{Factoring Groups into Subsets}.
Lecture Notes in Pure and Applied Mathematics, vol.~257.
Chapman \& Hall/CRC, Boca Raton, FL, 2009.
\href{https://doi.org/10.1201/9781420090475}
     {doi:10.1201/9781420090475}.

\end{thebibliography}
\end{document}